\documentclass[10pt,a4paper,twoside]{article}
\usepackage{epsfig,color,graphics,graphicx}
\usepackage{amssymb,amsbsy,amsmath,amsfonts,amssymb,amscd}
\usepackage{latexsym,euscript,exscale,epsfig}
\usepackage[english]{babel}
\newtheorem{theorem}{Theorem}

\newtheorem{corollary}[theorem]{Corollary}

\newtheorem{lemma}[theorem]{Lemma}

\newtheorem{remark}{Remark}
\newtheorem{proposition}[theorem]{Proposition}

\newenvironment{proof}[1][Proof]{\noindent\textbf{#1.} }{\ \rule{0.5em}{0.5em}}

\newcommand{\m}{{\mathbf m}}
\newcommand{\x}{{\mathbf x}}

\newcommand{\ba}{{\mathbf a}}
\newcommand{\bb}{{\mathbf b}}

\newcommand{\bt}{{\boldsymbol \tau}}
\newcommand{\bn}{{\boldsymbol \nu}}

\newcommand{\D}{\mathrm{div }}
\newcommand{\esssup}{\hbox{\rm ess}\sup }
\newcommand{\essinf}{\hbox{\rm ess}\inf }

\begin{document}

\title{Boundary stabilization and control of wave equations\\
by means of a general multiplier method.}
\author{
{\sc Pierre CORNILLEAU}\thanks{Universit\'e de Lyon, \'Ecole centrale de Lyon,
D.M.I., Institut Camille-Jordan (C.N.R.S. U.M.R. 5208),
36 avenue Guy-de-Collongue, 69134 \'Ecully cedex, France.} ,
{\sc Jean-Pierre LOH\'EAC}$^*$\thanks{Independent Moscow University,
Laboratoire J.-V. Poncelet (C.N.R.S. U.M.I. 2615),
Bol. Vlasyevsky Per. 11, 119002 Moscow, Russia.} .
}
\date{}

\maketitle



\medskip

\begin{abstract}
We describe a general multiplier method to obtain boundary
stabilization of the wave equation by means of a (linear or quasi-linear)
Neumann feedback. This also enables us to get Dirichlet boundary control of
the wave equation.
This method leads to new geometrical cases concerning the
"active" part of the boundary where the feedback (or control) is applied.
\par \noindent Due to mixed boundary conditions, the Neumann feedback case
generate singularities. Under a simple geometrical condition concerning the orientation of the
boundary, we obtain a stabilization result in linear or quasi-linear cases.
\end{abstract}

\section*{Introduction}
In this paper we are concerned with control and stabilization of the wave
equation in a multi-dimensional body $\Omega\subset \mathbb{R}^n$.

Stabilization is obtained using a feedback law given by some
part of the boundary of the spacial domain and some function defined on this
part. The problem can be written as follows
$$ \left\{ 
\begin{matrix}
u''-\Delta u=0\, \hfill
&\text{in } \Omega \times \mathbb{R}^*_+ \, ,\hfill \\
u=0\, \hfill
&\text{on } \partial \Omega_D \times \mathbb{R}^*_+ \, ,\hfill \\
\partial_\nu u=F \, \hfill
&\text{on } \partial \Omega_N \times \mathbb{R}^*_+ \, ,\hfill \\
u(0)=u_0 \, \hfill
&\text{in } \Omega \, ,\hfill \\
u'(0)=u_1 \, \hfill
&\text{in } \Omega \, ,\hfill
\end{matrix}
\right. $$
where we denote by $u' $, $u'' $, $\Delta u$ and $\partial_\nu u$ the first time-derivative
of $u$, the second time-derivative of the scalar function $u$, the standard Laplacian of $u$ and the normal outward derivative
of $u$ on $\partial \Omega $, respectively; $(\partial \Omega_D $, $\partial \Omega_N) $ is a partition of $\partial \Omega$ and $F$ is the feedback function which may depend on state $(u,u')$, position $x$ and time $t$.\\
Our purpose here is to choose the feedback law, that is to say the feedback
function $F$ and the ``active'' part of the boundary, $\partial \Omega_N $, so
that for every initial data, the energy function
$$E(t)=\frac{1}{2}\int_\Omega (|u'(t)|^2 +|\nabla u(t)|^2 )\, d\x \, ,$$
is decreasing with respect to time $t$, and vanishes as $t\longrightarrow
\infty $.\\
Formally, we can write the time-derivative of $E$ as follows
$$E'(t)=\int_{\partial \Omega_N } F u' \, d\sigma \, ,$$
and a sufficient condition for $E$ to be non-increasing would be: 
$\, F u'\le 0\, $ on $\partial \Omega_N $.\\
Thanks to the multiplier method introduced by L.F. Ho \cite{Ho} in the framework of Hilbert Uniqueness Method \cite{L}, it can be shown that the energy function is uniformly decreasing as time
$t$ tends to $\infty $ by choosing
$\x \mapsto \m (\x )=\x -\x_0$, where $\x_0$ is some given point in $\mathbb{R}^n$ and
$$\partial \Omega_N =\left\{ \x \in \partial \Omega \, / \, 
\m (\x ).\bn (\x )>0\, \right\} \, ,\quad
F=-\m .\bn \, u'\, ,$$
where $\bn $ is the normal unit vector pointing outward of $\Omega $.\\
This method has been performed by many authors, see for instance Komornik and
Zuazua \cite{KZ}, Komornik \cite{Ko} and the references therein. 
Here we extend the above result for rotated multipliers defined in \cite{O2}
and we follow  the analysis of singularities initiated by Grisvard
\cite{Gr, Gr2} and extended by Bey, Loh\'eac and Moussaoui \cite{BLM}.
This last work leads to results in case of higher dimensional domains with a 
non-empty boundary interface 
$\Gamma =\overline{\partial \Omega_N }\cap \overline{\partial \Omega_D }$ 
under an additional geometrical assumption concerning the orientation of the 
boundary.
\par Concerning the control problem, our goal is to find $v$ such that the 
solution of
\[\left\{ \begin{array}{l}
     u^{\prime \prime} - \Delta u = 0\\
     u = 0\\
     u = v\\
     u (0) = u_0\\
     u^{^{\prime}} (0) = u_1
   \end{array} \right. \left. \begin{array}{l}
     \text{in }  \Omega \times (0, T),\\
     \text{on }  \partial \Omega_D \times (0, T),\\
     \text{on }  \partial \Omega_N \times (0, T),\\
     \text{in }  \Omega,\\
     \text{in } \Omega,
   \end{array} \right. \]
reaches an equilibrium at $t=T$.
\par \noindent We here follow \cite{Ho}: in this work, Ho used the multiplier technique. His main purpose was to prove an inverse 
inequality for the linear wave equation implying its exact controllability. He 
introduced the so-called \textit{exit condition}: the control region must
contain a subset of the boundary where the scalar product between the outward 
normal and the vector pointing from some origin towards the normal is positive.
By varying the origin, a family of boundary controls satisfying the condition
is obtained.
\par \noindent In the last decades, micro-local techniques and geometric optics analysis allowed to find geometrical characterization of control and minimal control time in the exact controllability of waves. This condition has been introduced in \cite{BLR} under the name of \textit{Geometric Control Condition (GCC)}. It generalized the previous exit condition.\\
There is a certain balance: with {\it GCC}, control time is optimal but the observability constant is not explicit. With {\it exit condition}, time is not optimal but observability constants can be explicit, which is very useful in theoretical and numerical estimations.
\par \noindent In this paper we extend the family of multipliers recently introduced by Osses \cite{O2}.

\section{Notations and main results}
Let $\Omega $ be a bounded open connected set of $\mathbb{R}^{n}(n\geq 2)$ such that 
\begin{equation}\label{G}
\partial \Omega \text{ is }\mathcal{C}^{2} \text{ in the sense of Ne\v{c}as \cite{Ne}}.
\end{equation}

\noindent In the sequel, we denote by $\mathrm{I}$ the $n \times n$ identity matrix and by $\mathrm{A}^s$ the symmetric part of a matrix $\mathrm{A}$.
Let $\m$ be a $\mathrm{W}^{1,\infty}(\Omega)$ vector-field such that 
\begin{equation}\label{m}
\essinf_{\overline{\Omega}} \bigl( \D(\m) \bigr) >\esssup_{\overline{\Omega}} \bigl( \D(\m)-2 \lambda_\m \bigr) 
\end{equation}
where $\D $ is the usual divergence operator and $\lambda_\m(\x )$ is the smallest eigenvalue of the real symmetric matrix $\nabla \m (\x )^s$. Using Sobolev embedding, one may also assume that $\m \in \mathcal{C}(\overline \Omega).$ 

\begin{remark}
The set of all $\mathrm{W}^{1,\infty}(\Omega)$ vector-fields such that \eqref{m} holds is an open cone. If $\m$ belongs to this set, we denote 
$$c(\m)=\frac{1}{2} \Bigl( \essinf_{\overline{\Omega}} \bigl( \D(\m)\bigr) -\esssup_{\overline{\Omega}}\bigl( \D(\m)-2\lambda_\m \bigr) \Bigr) \, .$$ 
\end{remark}

\subsubsection*{Examples}

\begin{itemize}
\item An affine example is given by 
$$\m (\x )= (\mathrm{A_1} + \mathrm{A_2}) (\x - \x_0)\, ,$$ 
where $\mathrm{A_1}$ is a definite positive matrix, $\mathrm{A_2}$ a skew-symmetric matrix and $\x_0$ any point in $\mathbb{R}^n$.
\item A non linear example is 
$$\m (\x ) = (d\mathrm{I + A}) (\x - \x_0) + \mathbf{F} (\x)\ , $$
where $d > 0$, $\mathrm{A}$ is a skew-symmetric matrix, $\x_0$ any point in $\mathbb{R}^n$ and $\mathbf{F}$ is a $\mathrm{W}^{1,\infty}(\Omega)$ vector field such that 
$$\esssup_{\overline{\Omega}} \| \nabla \mathbf{F}^s \| <\frac{d}{n}\, ,$$ 
where $\|\cdot\|$ stands for the usual $2$-norm of matrices. 
\end{itemize}

\noindent We consider a partition $(\partial \Omega_N, \partial \Omega_D)$ of $\partial \Omega $ such that
\begin{equation}\label{R}
\left|
\begin{array}{c}
\Gamma =\overline{\partial \Omega }_{D}\cap \overline{\partial \Omega }_{N} \text{\ is a }\mathcal{C}^{3}\text{-manifold of dimension }n-2 \text{ such that } \m.\bn=0 \text{ on } \Gamma,\\
\exists \ \omega \text{ neighborhood of } \Gamma \text{ such that } \partial \Omega \cap \omega \ \text{\ is a }\mathcal{C}^{3}\text{-manifold
of dimension }n-1,\\
\qquad \mathcal{H}^{n-1}(\partial \Omega_{D})>0 \ (\mathcal{H}^{n-1} \text{ is the $(n-1)$-dimensional Hausdorff measure}).
\end{array}
\right.
\end{equation}
Furthermore, we assume
\begin{equation}\label{S1}
\partial \Omega _{N}\subset\{ \x\in \partial \Omega \, /\, \m (\x ).\bn (\x )
\geq 0\} \, ,\quad
\partial \Omega _{D}\subset\{ \x\in \partial \Omega \, /\, \m (\x ).\bn (\x )
\leq 0\} \, .
\end{equation}
This assumption clearly implies:\quad $\m.\bn=0 \text{ on } \Gamma$.

\subsection*{Boundary stabilization}

\noindent Let $g:\mathbb{R}\rightarrow\mathbb{R}$ be a measurable function such that
\begin{equation}\label{F1}
g \text{ is non-decreasing}\quad \text{and} \quad \exists k_{+}>0\, :\ |g(s)|\leq k_{+}|s| \, \text{ a.e.} \, .
\end{equation}
Let us now consider the following wave problem,
\begin{equation*}
(S)\quad \left\{ 
\begin{array}{l}
u^{\prime \prime }-\Delta u=0 \\ 
u=0 \\ 
\partial _{\nu }u=-\m.\bn \, g(u^{\prime }) \\ 
u(0)=u_{0} \\ 
u^{^{\prime }}(0)=u_{1}
\end{array}
\right. \left. 
\begin{array}{l}
\text{in }\Omega \times \mathbb{R}_{+}^{\ast }\, , \\ 
\text{on }\partial \Omega _{D}\times \mathbb{R}_{+}^{\ast }\, , \\ 
\text{on }\partial \Omega _{N}\times \mathbb{R}_{+}^{\ast }\, , \\ 
\text{in }\Omega \, , \\ 
\text{in }\Omega \, ,
\end{array}
\right.
\end{equation*}
where initial data satisfy
\begin{equation*}
(u_0 ,u_1 )\in \mathrm{H}_D^1 (\Omega )\times \mathrm{L}^2 (\Omega )
\end{equation*}
with $\mathrm{H}_D^1 (\Omega )=\{v\in \mathrm{H}^1 (\Omega )\, :\,
v=0\text{ on }\partial \Omega_D \}$.

\noindent Problem $(S)$ is well-posed in this space. Indeed, following Komornik \cite{Ko}, we define the non-linear operator $\mathcal{W}$ on 
$\mathrm{H}_D^1 (\Omega )\times \mathrm{L}^2 (\Omega )$ by 
$$ \begin{array}{l}
\mathcal{W} (u,v)=(-v,-\Delta u)\, ,\\
D(\mathcal{W} )=\{(u,v)\in \mathrm{H}_D^1 (\Omega )\times \mathrm{H}_D^1
(\Omega )\, /\, \Delta u\in \mathrm{L}^2 (\Omega )\ \text{ and } \ 
\partial_\nu u=-\m.\bn \, g(v)\, 
\text{ on }\partial \Omega_N \}\, ,
\end{array} $$
so that $(S)$ can be written as follows,
\begin{equation*}
\left\{ 
\begin{array}{l}
(u,v)' +\mathcal{W} (u,v)=0\, , \\ 
(u,v)(0)=(u_0 ,u_1 )\, .
\end{array}
\right.
\end{equation*}
It is classical that $\mathcal{W}$ is a maximal-monotone operator on 
$\mathrm{H}_D^1 (\Omega )\times \mathrm{L}^2 (\Omega )$ and that
$D(\mathcal{W} )$ is dense in 
$\mathrm{H}_D^1 (\Omega )\times \mathrm{L}^2 (\Omega )$ for the usual norm. 
Following Br\'{e}zis \cite{B}, we can deduce that for any initial data
$(u_0 ,v_0 )$ in $D(\mathcal{W} )$ there is a unique strong solution $(u,v)$
such that $u\in \mathrm{W^{1,\infty }(\mathbb{R};H_{D}^{1}(\Omega ))}$ and 
$\Delta u\in \mathrm{L^{\infty }(\mathbb{R}_{+};L^{2}(\Omega ))}$. Moreover, for 
two initial data, the corresponding solutions satisfy
\begin{equation*}
\forall t\geq 0\, ,\quad \Vert (u^{1}(t),v^{1}(t))-(u^{2}(t),v^{2}(t))\Vert
_{\mathrm{H}_D^1 (\Omega )\times \mathrm{L}^2 (\Omega )}\leq
\Vert (u^{1}_{0},v^{1}_{0})-(u^{2}_{0},v^{2}_{0})
\Vert _{\mathrm{H}_D^1 (\Omega )\times \mathrm{L}^2 (\Omega )}\, .
\end{equation*}
Using the density of $D(\mathcal{W} )$, one can extend the map
\begin{eqnarray*}
D(\mathcal{W} )
&\longrightarrow
&\mathrm{H}_D^1 (\Omega )\times \mathrm{L}^2 (\Omega )\\
(u_0, v_0) 
&\longmapsto 
&(u (t), v (t))
\end{eqnarray*}
to a strongly continuous semi-group of contractions $(S(t))_{t\geq 0}$ and define
for $(u_{0},u_{1})\in \mathrm{H}_D^1 (\Omega )\times \mathrm{L}^2 (\Omega )$
the weak solution $(u(t),u'(t))=S(t)(u_{0},u_{1})$
with the regularity $u\in \mathcal{C}(\mathbb{R}_{+};\mathrm{H}_{D}^{1}(\Omega ))\cap \mathcal{C}^{1}(\mathbb{R}_{+}; \mathrm{L}^{2}(\Omega ))$.
We hence define the energy function of solutions by
$$E(0)=\frac{1}{2}\int_{\Omega }(|u_{1}|^{2}+|\nabla u_{0}|^{2})\, d\x
\quad \hbox{and} \quad
E(t)=\frac{1}{2}\int_{\Omega }(|u^{\prime }(t)|^{2}+|\nabla u(t)|^{2})\, d\x \, ,\, \text{ for }t>0\, .$$
In order to get stabilization results, we need further assumptions
concerning the feedback function $g$
\begin{equation}\label{F2}
\exists p\geq 1 \, ,\ \exists k_{-}>0\, :\quad
|g(s)|\ge k_{-} \min \{ |s|,|s|^p \}\, ,\ \hbox{a.e.} \, ,
\end{equation}
and the additional geometric assumption
\begin{equation}\label{S2}
\m .\bt \le 0 \, \quad \text{on } \Gamma \, ,
\end{equation}
where $\bt (\x )$ is the normal unit vector pointing outward of $\partial
\Omega_N $ at a point $\x\in \Gamma $ when considering $\partial \Omega_N $
as a sub-manifold of $\partial \Omega $.

\begin{remark} 
It is not necessary to assume that 
$$\mathcal{H}^{n-1}(\{ \x \in \partial \Omega_{N}\, / \, \m (\x ).\bn( \x )>0
\})>0 $$
to get stabilization. In fact, our choices of $\m$ imply such properties 
(see examples in Section 5) whether the energy tends to zero.
\end{remark}

\noindent A main tool in this work is Rellich type relations \cite{Re}.
They lead to results of controllability and stabilization for the wave problem (see \cite{KZ} and \cite{Ho}). When the interface $\Gamma $ is not empty, the key-problem is
to show the existence of a decomposition of the solution in a regular and a singular parts (see \cite{Gr,KMR}) in any dimension. The
first results towards this direction are due to Moussaoui \cite{Mo}, and
Bey-Loh\'{e}ac-Moussaoui \cite{BLM}.\\
In this new case, our goal is to generalize those Rellich relations.
This will lead us to get a stabilization result about $(S)$ under \eqref{S1}, \eqref{S2}.
As well as in \cite{Ko}, we shall prove here two results of uniform boundary stabilization.

\subsubsection*{Exponential boundary stabilization}

We here consider the case when $p=1$ in \eqref{F2}. This is satisfied when
$g$ is linear,
$$\exists \alpha >0 \, : \qquad \forall s\in \mathbb{R} \, ,\quad g(s)=\alpha s\, .$$
In this case, the energy function is exponentially decreasing.
\begin{theorem}\label{T1}
Assume that conditions \eqref{G}, \eqref{m}, \eqref{R} and \eqref{S1} hold and that the feedback function $g$
satisfies \eqref{F1} and \eqref{F2} with $p=1$.\\
Then under the further geometrical assumption \eqref{S2}, there exist $C>0$
and $T>0$
such that for every initial data in $\mathrm{H}_D^1 (\Omega )\times \mathrm{L}^2 (\Omega )$, the energy of the
solution $u$ of $(S)$ satisfies
\begin{equation*}
\forall t>T\, ,\quad
E(t)\le E(0)\, \exp \Bigl( 1-\frac{t}{C} \Bigr) \, .
\end{equation*}
The above constants $C$ and $T$ do not depend on initial data.
\end{theorem}

\subsubsection*{Rational boundary stabilization}

We here consider the case $p>1$ and we get rational boundary stabilization.
\begin{theorem}\label{T2}
Assume that conditions  \eqref{G}, \eqref{m}, \eqref{R} and \eqref{S1} hold and that the feedback function $g$
satisfies \eqref{F1} and \eqref{F2} with $p>1$.\\
Then under the further geometrical assumption \eqref{S2}, there exist $C>0$
and $T>0$
such that for every initial data in $\mathrm{H}_D^1 (\Omega )\times \mathrm{L}^2 (\Omega )$, the energy of the
solution $u$ of $(S)$ satisfies
\begin{equation*}
\forall t>T\, ,\quad
E(t)\le C\, t^{2/(p-1)} \, .
\end{equation*}
where $C$ depends on the initial energy $E(0)$.
\end{theorem}

\begin{remark}
Taking advantage of the work by Banasiak and Roach \cite{BR} who generalized
Grisvard's results \cite{Gr} in the piecewise regular case, we will see that Theorems \ref{T1} and \ref{T2} remain true in the bi-dimensional case when
assumption \eqref{G} is replaced by following one,
\begin{equation}\label{G'}
\begin{matrix}
&\partial \Omega \text{ is a curvilinear polygon of class } {\cal C}^2 \, ,
\hfill \\
& \text{each component of } \partial \Omega \setminus \Gamma \text{ is a }
{\cal C}^2 \hbox{-manifold of dimension } 1\, ,\\
\end{matrix}
\end{equation}
and when condition \eqref{S2} is replaced by
\begin{equation}\label{S2'}
\begin{matrix}
\forall \x \in \Gamma \, ,\quad  0\le \varpi_\x \le \pi \quad \text{and} \quad
\text{if }
\varpi_\x =\pi\, ,\, \, \m(\x ).\bt(\x )\le 0 \, .
\end{matrix}
\end{equation}
where $\varpi_\x$ is the angle of the boundary at point $\x$.
\end{remark}

\subsection*{Boundary control problem}

Our problem consists in finding $T_0$ such that for each $T > T_0$ and for
every $(u_0, u_1) \in \mathrm{L}^2 (\Omega) \times \mathrm{H}^{- 1} (\Omega)$,
there exists $v \in \mathrm{L}^2 (\partial \Omega_N \times (0, T))$ in such a 
way that the solution of the wave equation
\[ (\Sigma) \quad \left\{ \begin{array}{l}
     u^{\prime \prime} - \Delta u = 0\\
     u = 0\\
     u = v\\
     u (0) = u_0 \\
     u^{^{\prime}} (0) = u_1
   \end{array} \right. \left. \begin{array}{l}
     \text{in }  \Omega \times (0, T)\, ,\\
     \text{on }  \partial \Omega_D \times (0, T)\, ,\\
     \text{on }  \partial \Omega_N \times (0, T)\, ,\\
     \text{in }  \Omega \, ,\\
     \text{in } \Omega \, .
   \end{array} \right. \]
satisfies
\begin{equation}
  u (T) = u' (T) = 0 \quad \text{ in } \Omega .
\end{equation}

\begin{theorem}\label{T5}
Assume that \eqref{G}, \eqref{m}, \eqref{R} and \eqref{S1} hold.\\
Then if $\displaystyle{ T > 2 \frac{\left. \| \m \right\|_{\infty}}{c(\m)}}$, for every initial data $(u_0, u_1 ) \in \mathrm{L}^2 (\Omega) \times \mathrm{H}^{- 1} (\Omega)$, there exists a control function $v \in \mathrm{L}^2 (\partial \Omega_N \times (0, T))$ such that the corresponding solution of $(\Sigma)$ satisfies final condition $(10)$.
\end{theorem}

\noindent Our paper is organized as follows.\\
In Section 2, we extend
Rellich relations (Theorems \ref{T3} and \ref{T4}) for elliptic problems with mixed boundary
conditions.\\
In Section 3, we apply these relations to prove some stabilization results with linear or
quasi-linear Neumann feedback (Theorems  \ref{T1} and \ref{T2}).\\
In Section 4, we extend some
observability and controllability results for the wave equation (Proposition \ref{obs} and Theorem \ref{T5}).\\
In Section 5, we detail affine examples in the case of a square domain.

\section{Rellich relation}

Here, we briefly extend Rellich relation obtained in \cite{BLM}, \cite{CLO} to our framework.

\subsection{A regular case}

We can easily build a Rellich relation corresponding to the above vector-field $\m$ when considered functions are smooth enough.

\begin{proposition}\label{RR}
Assume that $\Omega $ is a open set of $\mathbb{R}^{n}$ with
boundary of class $\mathcal{C}^{2}$ in the sense of Ne\v{c}as.
If $u$ belongs to $\mathrm{H}^2 (\Omega )$ then
\begin{equation*}
 2 \int_\Omega \Delta u \, \m. \nabla u\, d\x =
\int_\Omega (\D(\m)\mathrm{I} - 2 (\nabla \m)^s) (\nabla u, \nabla u)\, d\x
+ \int_{\partial \Omega} (2 \partial_{\nu} u \, \m. \nabla u 
-\m. \bn \, |\nabla u|^2)\, d\sigma \, .
\end{equation*}
\end{proposition}

\begin{proof}
Using Green-Riemann identity we get
\[ \int_{\Omega_{\varepsilon}} \Delta u \, \m. \nabla u\, d\x =
\int_{\partial \Omega_{\varepsilon}} \partial_{\nu} u \, \m. \nabla u\, d\sigma
-\int_{\Omega_{\varepsilon}} \nabla u. \nabla (\m. \nabla u)\, d\x \, . \]
So, observing that $\displaystyle{ \nabla u. \nabla (\m. \nabla u)
=\frac{1}{2} \m .\nabla (|\nabla u|^2 )+\nabla u.(\nabla \m )\nabla u 
=\frac{1}{2} \m .\nabla (|\nabla u|^2 )+(\nabla \m )^s (\nabla u, \nabla u)}$,
for smooth functions $u$, we get
\begin{equation*}
2 \int_{\Omega_{\varepsilon}} \Delta u \, \m. \nabla u\, d\x =
\int_{\partial \Omega_{\varepsilon}} 2 \partial_{\nu} u\, \m. \nabla u\, d\sigma
-2\int_{\Omega_{\varepsilon}}(\nabla \m)^s (\nabla u,\nabla u)\, d\x 
-\int_{\Omega_{\varepsilon}} \m. \nabla (| \nabla u|^2 )\, d\x \, .
\end{equation*}
With another use of Green-Riemann formula, we obtain the required formula
thanks to a classical approximation.
\end{proof}

We will now try to extend this result to the case of a less regular element $u$ when $\Omega $ is smooth enough.

\subsection{Bi-dimensional case}

We begin by the plane case: it is the simplest case from the point of view of singularity theory, and its understanding dates from Shamir \cite{Sh}. 
\begin{theorem}\label{T3}
Assume $n=2$. Under the geometrical conditions \eqref{G'} and \eqref{R}, let $\mathrm{u\in H^1(\Omega )}$ such that
\begin{equation*}
\Delta u \in \mathrm{L}^2 (\Omega)\, , 
\quad u_{/\partial \Omega_D } \in \mathrm{H}^{3/2}(\partial \Omega_D )\, , 
\quad \partial_{\nu} u_{/\partial \Omega_N } \in \mathrm{H}^{1/2} 
(\partial \Omega_N )\, .
\end{equation*}
Then $\, 2\partial _{\bn }u(\m.\nabla u)-(\m.\bn )|\nabla u|^{2}\in
\mathrm{L}^{1}(\partial \Omega )$
and there exist some coefficients $(c_{\x })_{\x \in \Gamma }$ such that
\begin{eqnarray*}
2 \int_{\Omega} \Delta u\, \m .\nabla u\, d\x 
&=&
\int_{\Omega }(\D (\m )I-2(\nabla \m )^s) (\nabla u, \nabla u)\, d\x
+\int_{\partial \Omega } (2 \partial_{\nu} u\, \m .\nabla u-\m .\bn \, 
|\nabla u|^2 )\, d\sigma \\
& &
+\frac{\pi }{4} \sum_{\x /\varpi_\x =\pi } c_{\x }^{2}(\m .\bt )(\x )\, .
\end{eqnarray*}
\end{theorem}

\begin{proof}
We follow the proof of Theorem 4 in \cite{CLO} to get this result.
\end{proof}

\begin{remark}
As in Theorem 4 of \cite{CLO}, the assumption $\mathcal{H}^{1}(\partial \Omega_{D})>0$ is not necessary in the above proof.
\end{remark}

\subsection{General case}

We now state the result in higher dimension.
\begin{theorem}\label{T4}
Assume $n\geq 3$. Under geometrical conditions \eqref{G} and \eqref{R}, let
$u\in \mathrm{H}^1(\Omega )$ such that 
\begin{equation*}
\Delta u \in \mathrm{L}^2 (\Omega)\, , 
\quad u_{/\partial \Omega_D } \in \mathrm{H}^{3/2}(\partial \Omega_D )\, , 
\quad \partial_{\nu} u_{/\partial \Omega_N } \in \mathrm{H}^{1/2} (\partial
\Omega_N )\, .
\end{equation*}
Then 
$\, 2\partial _{\bn }u(\m.\nabla u)-(\m.\bn )|\nabla u|^{2}\in \mathrm{L}^{1}
(\partial \Omega )$ and there exists $\zeta \in \mathrm{H}^{1/2}(\Gamma )$ 
such that
\begin{eqnarray*}
2\int_{\Omega } \Delta u \, \m .\nabla u\, d\x
&=&
\int_{\Omega } (\D (\m )I-2(\nabla \m )^s)(\nabla u, \nabla u) \, d\x
+\int_{\partial \Omega } (2 \partial_{\nu} u\, \m. \nabla u-\m .\bn \,
|\nabla u|^2 )\, d\sigma \\
& &
+\int_{\Gamma} \m. \bt \, |\zeta |^2 \, d\gamma \, .
\end{eqnarray*}
\end{theorem}

\begin{proof}
We exactly follow the proof of Theorem 5 in \cite{CLO} to get this result.
\end{proof}

\section{Linear and quasi-linear stabilization}

We begin by writing the following consequence of results of Section 2.
\begin{corollary}\label{C}
Assume that $t\mapsto(u(t),u'(t))$ is a strong solution of $(S)$ and that geometrical additional assumption \eqref{S2} if $n\geq3$ (or \eqref{S2'} if $n=2$) holds, then, for every time $t$, $u(t)$ satisfies
\begin{equation*}
2\int_{\Omega }\Delta u\, \m .\nabla u\, d\x \leq
\int_{\Omega} (\D (\m )I-2(\nabla \m )^s)(\nabla u,\nabla u ) \, d\x 
+\int_{\partial \Omega } (2 \partial_{\nu} u\, \m .\nabla u-\m .\bn \,
|\nabla u|^2)\, d\sigma \, .
\end{equation*}
\end{corollary}
\begin{proof}
 Indeed, under theses hypotheses, for each time $t$, $(u(t),u'(t))\in D(\mathcal{W})$ so that $u(t)$ satisfies  hypotheses of Theorems \ref{T3} or \ref{T4}. The result follows immediately from \eqref{S2} or \eqref{S2'}.
\end{proof}

The main tool in the proof of Theorems \ref{T1}, \ref{T2} is the following result (see proof in \cite{Ko}) which will be applied with $\displaystyle{\alpha=\frac{p-1}{2}}$.

\begin{proposition}\label{K}
Let $E:\mathbb{R}_{+}\rightarrow \mathbb{R}_{+}$ be a non-increasing function such that there exist $\alpha \geq 0$ and $C>0$ which fulfill
\begin{equation*}
\forall t\geq 0\text{, }\int_{t}^{\infty }E^{\alpha +1}(s)ds\leq CE(t).
\end{equation*}
Then, \ setting $T=CE^{\alpha }(0)$, one gets
\begin{eqnarray*}
    \text{if } \alpha=0\, ,\quad  &\forall t \geq T,& \, E (t) \leq E (0) \exp
    \left( 1 - \frac{t}{T}\right) \, ,\\
    \text{if } \alpha>0\, ,\quad  &\forall t \geq T,& \, E (t) \leq E (0) \left(
    \frac{T + \alpha T}{T + \alpha t} \right)^{1/\alpha} \, .
  \end{eqnarray*}
\end{proposition}

As usual in this context, we will perform the multiplier method to a suitable $\m$.\\
Putting $Mu=2\m.\nabla u+a u$ with $a$ a constant to be defined later, we prove the following result.

\begin{lemma}\label{IPP}
For any $0\leq S<T<\infty $, the following inequality holds
\begin{eqnarray*}
\displaystyle{ \int_S^T E^{\frac{p - 1}{2}} \int_{\Omega} \Bigl( (\D(\m )-a) 
(u')^2 +\bigl( (a - \D (\m ))I+2(\nabla \m )^s \bigr) (\nabla u,\nabla u)\Bigr)
\, d\x \, dt } \\
\le \displaystyle{- \Bigl[ E^{\frac{p - 1}{2}} \int_{\Omega} u^{\prime} Mu\, d\x 
\Bigr]_S^T 
+\frac{p - 1}{2} \int_S^T E^{\frac{p - 3}{2}} E^{\prime} \int_{\Omega} u^{\prime}
Mu\, d\x \, dt} \\
\displaystyle{+ \int_S^T E^{\frac{p - 1}{2}} \int_{\partial \Omega_N} \m .\bn \,
\bigl( (u^{\prime})^2 -|\nabla u|^2 -g(u^{\prime}) Mu\bigr) \, d\sigma \, dt
\, .}
\end{eqnarray*}
\end{lemma}
  
\begin{proof}
We here follow \cite{CR}.\\
We Use the fact that $u$ is solution of $(S)$ and 
we observe that $u^{\prime
\prime }Mu=(u^{\prime }Mu)^{\prime }-u^{\prime }Mu^{\prime }$. Then an
integration by parts gives
    \begin{eqnarray*}
     0 &=&\int_S^T E^{\frac{p - 1}{2}} \int_{\Omega} (u^{\prime \prime} -
      \Delta u)\, Mu \, d\x \, dt\\
      &=&\Bigl[ E^{\frac{p - 1}{2}} \int_{\Omega} u^{\prime} Mu\, d\x
      \Bigr]_S^T - \frac{p - 1}{2} \int_S^T E^{\frac{p - 3}{2}} E^{\prime}
      \int_{\Omega} u^{\prime} Mu\, d\x \, dt
      - \int_S^T E^{\frac{p - 1}{2}} \int_{\Omega} (u^{\prime}
      Mu^{\prime} + \Delta u\, Mu) \, d\x \, dt\, .
    \end{eqnarray*}
Corollary \ref{C} now gives
$$\int_{\Omega} \Delta u\, Mu\, d\x \leq  a \int_{\Omega}  \Delta u\, u \, d\x 
+\int_{\Omega} (\D (\m )I-2(\nabla \m)^s) (\nabla u, \nabla u)\, d\x
+\int_{\partial \Omega} (2 \partial_{\nu} u\, \m .\nabla u-\m .\bn \,
|\nabla u|^2)\, d\sigma \, .$$
Consequently, Green-Riemann formula leads to
\[ \int_{\Omega} \Delta u\, Mu\, d\x \leq \int_{\Omega} ((\D (\m )-a)I-2(\nabla
\m)^s) (\nabla u, \nabla u)\, d\x +
\int_{\partial \Omega} (\partial_{\nu} u\, Mu-\m .\bn \, |\nabla u|^2)\, d\sigma
\, . \]
Using boundary conditions and the fact that $\nabla u =\partial_{\nu} u\,
\bn$ on $\partial \Omega_D$, we then get
\[ \int_{\Omega} \Delta uMu\, d\x \leq \int_{\Omega} ((\D (\m) -
       a) I - 2 (\nabla \m)^s) (\nabla u, \nabla u) \, d\x - \int_{\partial
       \Omega_N} \m. \bn \, (g (u^{\prime})\, Mu + | \nabla u|^2)\, d\sigma \, . \]
On the other hand, another use of Green formula gives us
    \[ \int_{\Omega} u^{\prime} Mu^{\prime} \, d\x = \int_{\Omega} (a -
       \D (\m)) (u^{\prime})^2 \, d\x + \int_{\partial \Omega_N} \m. \bn \,
       |u^{\prime} |^2 \, d\sigma \, . \]
We complete the proof by summing up above estimates.
\end{proof}

\noindent Let us now prove Theorems \ref{T1} and \ref{T2}.

\begin{proof}
Following \cite{Ko} and \cite{CR}, we will prove the estimates for $(u_{0},u_{1})\in D(\mathcal{W})$ which will be sufficient thanks to a density argument. 

\noindent Using Lemma \ref{IPP}, we have to find $a$ such that $\D (\m) - a$ and $(a - \D (\m)) I + 2 (\nabla \m)^s$ are uniformly minorized on $\Omega$, that is, almost everywhere on $\Omega$
\begin{equation}\label{Paramètre c}
\left\{ \begin{array}{l}
     \D(\m)-a\geq c\, ,\\
     2\lambda_\m+(a-\D(\m))\geq c\, ,
   \end{array} \right.
\end{equation}
for some positive constant $c$.
The latter condition is then equivalent to find $a$ which fulfills 
$$\essinf_{\bar{\Omega}}\bigl( \D(\m )\bigr) >a>
\esssup_{\bar{\Omega}}\bigl( \D(\m)-2\lambda_\m \bigr) \, ,$$
and its existence is now garanted by \eqref{m}. Moreover, it is straightforward to see that the greatest value of $c$ such that \eqref{Paramètre c} holds is 
$$\frac{1}{2} \Bigl( \essinf_{\bar{\Omega}}\bigl( \D(\m)\bigr)
-\esssup_{\bar{\Omega}} \bigl( \D(\m)-2\lambda_\m \bigr) \Bigr) =c(\m )\, ,$$
and obtained for 
$\displaystyle{a=a_0:=\frac{1}{2} \Bigl(\essinf_{\bar{\Omega}} \bigl( \D(\m)\bigr) +\esssup_{\bar{\Omega}} \bigl( \D(\m)-2\lambda_\m\bigr) \Bigr) }$.\\
With this value $a_0$, we apply Lemma \ref{IPP} and get
\begin{eqnarray*}
2c(\m)\int_{S}^{T}E^{\frac{p+1}{2}}\, dt &\leq &-\Bigl[ E^{\frac{p-1}{2}%
}\int_{\Omega }u^{\prime } Mu\, d\x\Bigr]_{S}^{T}+\frac{p-1}{2}\int_{S}^{T}E^{
\frac{p-3}{2}}E^{\prime } \int_{\Omega }u^{\prime } Mu\, d\x \, dt \\
&&+\int_{S}^{T}E^{\frac{p-1}{2}} \int_{\partial \Omega _{N}} \m .\bn \,
\bigl( (u^{\prime })^{2}-|\nabla u|^{2}-g(u^{\prime })Mu \bigr)\, d\sigma \, dt\, .
\end{eqnarray*}
Young and Poincar\'{e} inequality gives
$$\Bigl| \int_{\Omega }u^{\prime }Mu\, d\x \Bigr| \leq CE(t)\, .$$ 
It follows then
  \begin{eqnarray*}
    2 c(\m) \int_S^T E^{\frac{p + 1}{2}} \, dt & \leq & C (E^{\frac{p + 1}{2}} (T) +
    E^{\frac{p + 1}{2}} (S)) + C \int_S^T E^{\frac{p - 1}{2}} E^{\prime} \, dt\\
    &  & + \int_S^T E^{\frac{p - 1}{2}} \int_{\partial \Omega_N} \m .
    \bn \, \bigl( (u^{\prime})^2 - | \nabla u|^2 - g (u^{\prime}) Mu\bigr)\, d\sigma
    \, dt \, .
  \end{eqnarray*}
  Let $d\sigma_{\m} =\m.\bn \, d\sigma $. If we observe that $\displaystyle{E^{\prime
}(t)=-\int_{\partial \Omega _{N}}g(u^{\prime })u^{\prime }\, d\sigma _{\m}\leq 0}$,
we get, for a constant $C>0$ independent of $E(0)$ if $p=1$,
  \[ 2 c(\m) \int_S^T E^{\frac{p + 1}{2}} \, dt \leq CE (S) + \int_S^T E^{\frac{p -
     1}{2}} \int_{\partial \Omega_N} \bigl( (u^{\prime})^2 - | \nabla u|^2 - g
     (u^{\prime}) Mu\bigr) \, d\sigma_\m \, dt. \]
 Using the definition of $Mu$ and Young inequality, we get for any $\varepsilon_0 > 0$
\[ 2 c(\m) \int_S^T E^{\frac{p + 1}{2}} \, dt \leq CE (S) + \int_S^T
E^{\frac{p - 1}{2}} \int_{\partial \Omega_N} \bigl( (u^{\prime})^2(1+\| \m\|_{\infty}^2) + \frac{\alpha^2}{4 \varepsilon_0} g (u^{\prime})^2 +\varepsilon_0 u^2\bigr) \, d\sigma_\m \, dt \, .\]
Now, using Poincar\'{e} inequality, we can choose $\varepsilon_0 >0$ such
that
  \[ \varepsilon_0 \int_{\partial \Omega_N} u^2\, d\sigma_\m \leq \frac{c(\m)}{2}
     \int_{\Omega} | \nabla u|^2 \, d\x \leq c(\m)E \, .\]
So we conclude 
\begin{equation*}
c(\m)\int_{S}^{T}E^{\frac{p+1}{2}}\, dt\leq CE(S)+C\int_{S}^{T}E^{\frac{p-1}{2}%
}\int_{\partial \Omega _{N}}\bigl( (u^{\prime })^{2}+g(u^{\prime
})^{2}\bigr)\, d\sigma _{\m} \, dt\, .
\end{equation*}
We split $\partial \Omega_{N}$ to bound the last term of this estimate
 $$\partial \Omega _{N}^{1}=\{\x\in
\partial \Omega _{N}\, /\, |u^{\prime }(\x)|>1\} \, , \quad \partial \Omega
_{N}^{2}=\{\x\in \partial \Omega _{N}\, /\, |u^{\prime }(\x)|\leq 1 \} \, .$$
Using \eqref{F1} and \eqref{F2}, we get
\begin{equation*}
\int_{S}^{T}E^{\frac{p-1}{2}} \int_{\partial \Omega_{N}^{1}}
\bigl( (u^{\prime })^{2}+g(u^{\prime })^{2}\bigr)\, d\sigma_{\m} \, dt \, \leq
\, C\int_{S}^{T}E^{\frac{p-1}{2}} \int_{\partial \Omega _{N}}u^{\prime}g(u^{\prime })\, d\sigma_{\m} \, dt \, \leq \, CE(S)\, ,
\end{equation*}
where $C$ depends on $E(0)$ if $p=1$.

\noindent On the other hand, using \eqref{F1}, \eqref{F2}, Jensen inequality and boundedness of $\m$, one successively obtains
\begin{equation*}
\int_{\partial \Omega _{N}^{2}}\bigl( (u^{\prime })^{2}+g(u^{\prime })^{2}\bigr) \, d\sigma_{m} \leq C\int_{\partial \Omega _{N}^{2}}(u^{\prime }g(u^{\prime
}))^{2/(p+1)}\, d\sigma_{\m} \leq C\Bigl( \int_{\partial \Omega _{N}^{2}}u^{\prime }g(u^{\prime})\, d\sigma_{\m}\Bigr)^{\frac{2}{p+1}} \leq C(-E^{\prime })^{\frac{2}{p+1}}\, .
\end{equation*}
Hence, using Young inequality again, we get for every $\varepsilon >0$
\begin{equation*}
\int_{S}^{T}E^{\frac{p-1}{2}} \int_{\partial \Omega _{N}^{2}}\bigl( (u^{\prime })^{2}+g(u^{\prime })^{2}\bigr)\, d\sigma_{\m} \, dt \leq
\int_{S}^{T}(\varepsilon E^{\frac{p+1}{2}}-C(\varepsilon )E^{\prime })\, dt
\leq \varepsilon \int_{S}^{T}E^{\frac{p+1}{2}}\, dt+C(\varepsilon )E(S)\, .
\end{equation*}
Finally we get, for some $C(\varepsilon )$ and $C$ independent of $E(0)$
if $p=1$
\begin{equation*}
c(\m)\int_{S}^{T}E^{\frac{p+1}{2}}\, dt \leq C(\varepsilon )E(S)+\varepsilon
C\int_{S}^{T}E^{\frac{p+1}{2}}\, dt \, .
\end{equation*}
Choosing now $\displaystyle{\varepsilon C \leq \frac{c(\m)}{2}}$, one obtains
$$
c(\m)\int_{S}^{T}E^{\frac{p+1}{2}}\, dt\leq CE(S)\, ,
$$
and Theorems can be deduced from Lemma \ref{K}.
\end{proof}

\begin{remark}
As stated before, we can replace $\m$ by $\lambda \m$ for any positive
 $\lambda$. One can wonder what happens to the speed of stabilization
 $\theta=\frac{c(\m)}{C}$ found in Theorem 3. In fact, a careful estimation
 of all terms shows that one can obtain
 \[ C = k_- + k_+ \lambda^2 + k_+ \frac{a_0^2}{4} (1 +
     C_P) C_{Tr} \lambda^3 \, , \]
where $C_P$ denotes the Poincar\'e constant and $C_{Tr}$ the norm of
the trace application $Tr: \mathrm{H}^1 (\Omega) \rightarrow \mathrm{L}^2 (\partial \Omega)$. The speed found in our proof is consequently
 $$\theta = c(\m) \left( \frac{k_-}{\lambda} + k_+
\lambda + k_+ \frac{a_0^2}{4} (1 + C_P) C_{Tr} \lambda^2
\right)^{- 1}\, .$$ 
It can be shown that $\theta$ reaches a maximum at some point 
$$ \lambda_0 \in \biggl[ \min \Bigl( \Bigl( \frac{k_-}{k_+ a_0^2 (1 + C_P) C_{Tr}}\Bigr) ^{1/3}, \frac{2}{a_0^2 (1 + C_P) C_{Tr}}
  \Bigl), \sqrt{\frac{k_-}{k_+}} \ \biggr] \, .$$
Besides, $\theta$ tends to $0$  when $\lambda \rightarrow 0$ or $\infty$. 
\end{remark}

\begin{remark}
  In fact, one can replace the feedback law $\m. \bn \, g (u')$ by a more general
  one $g (\x, u')$ provided that, for some constant $c>1$,
  \begin{eqnarray*}
   c^{-1} (\m. \bn)^{\frac{1}{p}} |s|^{\frac{1}{2} + \frac{1}{p}} \leq  & | g (\x, s)| & \leq c (\m. \bn)^{\frac{1}{p}} |s|^{\frac{1}{2} + \frac{1}{p}} \quad \text{for a.e. } \x \in \partial \Omega_N\text{ and } |s| \leqslant 1\, ,\\
     c^{-1} (\m. \bn) |s| \leq & | g (\x, s)| & \leq c (\m. \bn) |s| \qquad \quad \ \text{for a.e. } \x \in \partial \Omega_N \text{ and } |s|\geqslant 1\, . 
  \end{eqnarray*}
The details are left to the reader but the previous proof works also in
this case.
\end{remark}

\section{Observability and controllability results}

It is well-known that micro-local techniques \cite{BLR} characterize all partitions of
the boundary such that this result holds, but constants are not
explicit. Thus, using this new choice of multiplier, we will enlarge the set
of geometric examples with explicit knowledge of constants. We here follow \cite{O2}.

\subsection{Preliminary settings}

Following HUM method \cite{L}, controlabillity of problem $(\Sigma )$ is
equivalent to observability of its adjoint problem.
the solution of the control problem is equivalent to
studying the observability properties of the adjoint problem. For each pair
of initial conditions $(\varphi_0, \varphi_1) \in {\mathrm H}_0^1 (\Omega )
\times {\mathrm L}^2 (\Omega )$, let us consider the solution $\varphi $ of 
the following wave problem,
\[ (\Sigma' ) \quad \left\{ \begin{array}{l}
     \varphi^{\prime \prime} - \triangle \varphi = 0\\
     \varphi = 0\\
     \varphi (0) = \varphi_0\\
     \varphi^{\prime} (0) = \varphi_1
   \end{array} \right. \left. \begin{array}{l}
     \text{in }  \Omega \times (0, T)\, ,\\
     \text{on }  \partial \Omega \times (0, T)\, ,\\
     \text{in }  \Omega\, ,\\
     \text{in }  \Omega \, .
   \end{array} \right. \]
Observability of $(\Sigma' )$ is equivalent to the existence of a constant 
$C < \infty$ independent of $(\varphi_0, \varphi_1)$ such that
\begin{equation*}
 E_0 = \frac{1}{2} \int_{\Omega} \bigl( |\varphi_1 |^2 +|\nabla \varphi_0 |^2 
\bigr)\, d\x \leqslant C \int_{\partial \Omega_N \times (0, T)} | \partial_{\nu}
  \varphi |^2 d\sigma \, dt\, .
\end{equation*}
Let us define the operator $\mathcal{W}_0$ on 
$\mathrm{H}_0^1 (\Omega) \times \mathrm{L}^2 (\Omega )$ by 
$$\begin{matrix}
\mathcal{W}_0 (\varphi, \psi) = (- \psi, - \Delta
\varphi) \, ,\hfill \\
D(\mathcal{W}_0 )=\{(\varphi, \psi) \in \mathrm{H}_0^1 (\Omega) \times
\mathrm{H}_0^1 (\Omega) \, /\, \Delta \varphi \in \mathrm{L}^2 (\Omega) \}\, ,
\hfill 
\end{matrix} $$
so that $(\Sigma')$ can be written as follows,
\[ \left\{ \begin{array}{c}
     (\varphi, \psi)^{\prime} + \mathcal{W}_0 (\varphi, \psi) = 0 \, ,\\
     (\varphi, \psi) (0) = (\varphi_0, \varphi_1) \, .
   \end{array} \right. \]
\begin{remark}
If $(\varphi, \psi) \in D(\mathcal{W}_0 )$, $\varphi $ is the solution of 
some Dirichlet Laplace problem and hence regular (that is \ $\varphi \in
  \mathrm{H}^2 (\Omega)$).
\end{remark}

\noindent $\mathcal{W}_0 $ is a maximal-monotone operator on 
$\mathrm{H}_0^1 (\Omega) \times \mathrm{L}^2 (\Omega )$ and $D (\mathcal{W}_0 )$
is dense in $\mathrm{H}_0^1 (\Omega) \times \mathrm{L}^2 (\Omega )$
for the usual norm. Using Hille-Yosida Theorem, it generates a unitary 
semi-group on $\mathrm{H}_0^1 (\Omega) \times \mathrm{L}^2 (\Omega )$,
we denote its value applied at $(\varphi_0 ,\varphi_1 )$ at time $t$ by 
$(\varphi (t), \varphi' (t))$. 

\noindent As a consequence, we get conservation of energy.
\begin{proposition}
 If $t \geqslant 0$ and $\varphi$ is a weak solution of $(\Sigma')$, then
 \[  E (t) = \frac{1}{2} \int_{\Omega} (| \varphi^{\prime} (t) |^2 + |
       \nabla \varphi (t) |^2) \, d\x = E_0 \, . \]
\end{proposition}

\noindent A weak solution of $(\Sigma')$ hence belongs to
$\mathcal{C} ( \mathbb{R}_+ ; \mathrm{H}_D^1 (\Omega)) \cap \mathcal{C}^1 (
\mathbb{R}_+ ; \mathrm{L}^2 (\Omega))$.\\ A solution with 
$(\varphi_0, \varphi_1) \in D (\mathcal{W}_0 )$ is called a strong solution
and satisfies
$(\varphi, \varphi') \in \mathcal{C} ( \mathbb{R}_+ ;D(\mathcal{W}_0 ))$.

\subsection{Inverse inequality and exact controllability}

We keep similar notations as in Section 2: 
$\displaystyle{a_0=\frac{1}{2} \Bigl(\essinf_{\overline{\Omega}}\bigl( \D (\m )\bigr)
+\esssup_{\overline{\Omega}}\bigl( \D (\m )-2 \lambda_\m \bigr) \Bigr) }$.

\begin{proposition}\label{obs}
If $\displaystyle{T > 2 \frac{ \| \m \|_\infty }{c(\m )}}$, for each weak
solution $\varphi$ of $(\Sigma')$, the following inequality holds
\[ E_0 \leqslant \frac{\displaystyle{ \esssup_{\partial \Omega_N} }
|\m. \bn |}{ 2\bigl( c(\m ) T - 2 \| \m \|_\infty \bigr) }
\int_{\partial \Omega_N \times (0, T)} |\partial_\nu \varphi |^2 \, d\sigma \,
dt\, . \]
\end{proposition}

\begin{remark}
In the case $\m (x) = (d \mathrm{I + A}) (\x - \x_0)$ with $\mathrm{A}$ skew-symmetric matrix, we recover classical results (see \cite{KZ},\cite{O2}).
\end{remark}

\begin{proof}
Let $(\varphi_0, \varphi_1) \in D (\mathcal{W}_0 )$. We use again
$M \varphi =2 \m .\nabla \varphi + a_0 \varphi$. 
Using the fact that $\varphi $ is solution of $(\Sigma' )$ and observing that 
$\varphi^{\prime \prime } M\varphi =(\varphi' M \varphi )' -\varphi' M\varphi' $,
we get
$$0=\int_0^T \int_\Omega (- \varphi^{\prime \prime} +\Delta \varphi )M\varphi
\, d\x \, dt =
-\Bigl[ \int_\Omega \varphi' M\varphi \, d\x \Bigr]_0^T
+\int_0^T \int_\Omega (\varphi' M \varphi' + \Delta \varphi M\varphi )
\, d\x \, dt \, .$$
As well as in the proof of Theorems \ref{T1} and \ref{T2}, one uses Green-Riemann formula and Proposition \ref{RR} to get
\[  \int_{\Omega} \Delta \varphi M \varphi \, d\x = \int_{\Omega}
        ((\D (\m) - a_0) I - 2 (\nabla \m)^s) (\nabla \varphi,
        \nabla \varphi) \, d\x + \int_{\partial \Omega} (\partial_{\nu} \varphi
        M \varphi - \m. \bn \, | \nabla \varphi |^2)\,  d\sigma \, . \]
Dirichlet boundary conditions lead to
  \[ \int_{\Omega} \Delta \varphi M \varphi \, d\x =
           \int_{\Omega} ((\D (\m) - a_0) I - 2 (\nabla \m)^s)
           (\nabla \varphi, \nabla \varphi) \, d\x + \int_{\partial \Omega_N}
           \m. \bn \, |\partial_{\nu} \varphi |^2 \, d\sigma \, . \]
On the other hand, another use of Green formula gives us
  \[ \int_{\Omega} \varphi' M \varphi^{\prime} \, d\x = \int_{\Omega} (a_0
     - \D (\m)) | \varphi^{\prime} |^2 \, d\x \, , \]
so, we finally get, using the same minoration as in proof of Theorems
\ref{T1} and \ref{T2}
  \begin{equation}\label{EE}
    c(\m) \int_0^T \int_{\Omega} | \varphi^{\prime} |^2 + | \nabla \varphi
    |^2 \, d\x \, dt \leqslant - \Bigl[ \int_{\Omega} \varphi^{\prime} M \varphi \, d\x
    \Bigr]_0^T + \int_{\partial \Omega_N} \m. \bn \, | \partial_{\nu}
    \varphi |^2 \, d\sigma \, .
  \end{equation}
 Using the conservation of the energy, the left hand side in \eqref{EE} is $2 c T E_0$. It only
  remains to estimate the term $\displaystyle{\Bigl[ - \int_{\Omega} \varphi^{\prime} M
  \varphi \, d\x \Bigr]_0^T}$ to end the proof.
  
\noindent Let us fix a time $t \in \{ 0, T\} $.
Cauchy-Schwarz inequality leads to
  \[  - \int_{\Omega} \varphi^{\prime} M \varphi \, d\x \leqslant \Bigl(
           \int_{\Omega} |_{} \varphi^{\prime} |^2 \Bigr)^{1/2}
           \Bigl( \int_{\Omega} |_{} M \varphi |^2 \Bigr)^{1/2} \]
Denoting by $\| . \|$ the $\mathrm{L}^2 (\Omega)$-norm, we get the
following splitting
  \[ \| M \varphi \|^2 = \|2 m. \nabla \varphi
     \|^2 + a_0^2 \| \varphi \|^2 + 4 a_0 \int_\Omega
     \varphi \, \m .\nabla \varphi \, d\x \, . \]
 Green-Riemann formula and Dirichlet boundary conditions give
  \[ \int_\Omega
     \varphi \, \m .\nabla \varphi \, d\x = - \frac{1}{2}
     \int_{\Omega} \D (\m) | \varphi |^2 \, d\x \, ,\]
  and since $a_0 - 2 \D (\m) \leqslant a_0 -
  \D (\m) \leqslant - c(\m) \text{ a.e.}$, we finally get that $\| M \varphi \| \leqslant 2 \| \m \|_{\infty} \| \nabla \varphi \|$.
  
\noindent Consequently, with Young inequality, we get the following estimate
  \[ - \int_{\Omega} \varphi^{\prime} M \varphi \, d\x \leqslant 2 \| \m
     \|_{\infty} \Bigl( \int_{\Omega} |_{} \varphi^{\prime} |^2
     \Bigr)^{1/2} \Bigl( \int_{\Omega} | \nabla \varphi |^2
     \Bigr)^{1/2} \leqslant 2 \| \m \|_{\infty} E_0 \, . \]
So \eqref{EE} becomes
\[ 2 \bigl( c(\m )T-2\| \m \|_\infty \bigr) E_0\leqslant
\int_{\partial \Omega_N} \m. \bn \, |\partial_\nu \varphi |^2 \, d\sigma \, , \]
which ends the proof of Proposition \ref{obs}, using the density of the domain.
\end{proof}

\noindent Now we can deduce our exact controllability result (Theorem \ref{T5})
from Proposition \ref{obs} following HUM method (see \cite{L}, Chapter IV). 

\section{Example}
Let us consider here the case of a square domain $\Omega=(0,1)^2 $ with the
following affine multiplier,
\begin{equation}
 \m(\x)=\left( \begin{array}{c}
     \cot \theta_1\\
     1
   \end{array} \begin{array}{l}
     -1\\
     \cot \theta_2
   \end{array} \right)(\x-\x_0)
\end{equation}
where $\theta_1 $ and $\theta_2 $ belong to
$\displaystyle{ \left( 0 ,\frac{\pi }{2} \right) }$.\\ 
We will discuss the dependence of $\partial \Omega_N $ and $\partial \Omega_D $
on $\x_0 $.\\
First let us consider one edge $[\ba \bb ]$ of $\Omega $ with its normal unit
vector $\bn $. One can easily see that
$$\m (\x ).\bn (\x )=\frac{1}{\sin \theta } (\x -\x_0 ).\bn_\theta \, ,$$
where $\theta =\theta_1 $ (resp. $\theta_2 $) if $[\ba \bb ] \subset [0,1]
\times \{ 0,1\}$ (resp. $\{ 0,1\} \times [0,1] $) and $\bn_\theta $ is deduced
from $\bn $ by rotation of angle $-\theta $.
Without any restriction, we suppose
\quad $\ba .\bn_\theta < \bb .\bn_\theta $.\\
Then there exists an interface point along $[\ba \bb ] $ if and only if
$\x_0 $ belongs to the belt
$${\cal B}_\theta =\{ \x \in \mathbb{R}^2 \, /\, 
\ba .\bn_\theta < \x .\bn_\theta < \bb .\bn_\theta \} \, .$$
In this case, at this interface point $\x_1 $, we get with similar notations,
\quad
$$ \m (\x_1 ).\bt (\x_1 )=\frac{1}{\sin \theta } (\x_1 -\x_0 ).\bt_\theta \, .$$
Then additional geometric assumption \eqref{S2} is not satisfied if $\x_0 $
belongs to half-belt ${\cal B}_\theta^+ $ (see Fig. 1).
\begin{figure}[htb]
\centerline{
\begin{picture}(0,0)
\includegraphics{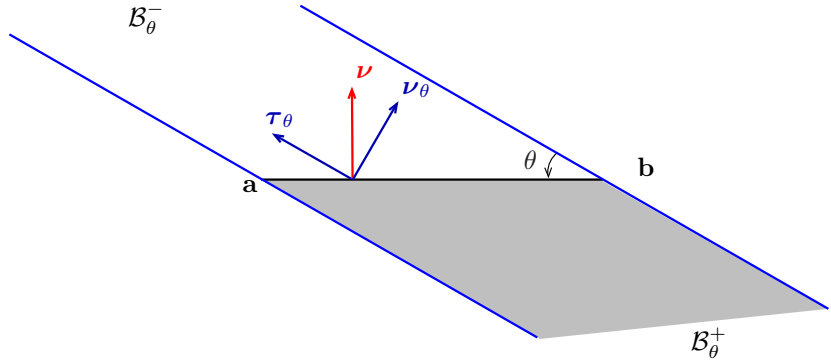}
\end{picture}
\setlength{\unitlength}{3729sp}
\begin{picture}(5410,2463)(873,-2439)
\put(4850,-1150){\makebox(0,0)[lb]{\smash{
{\color[rgb]{0,0,0}$\bb $}
}}}
\put(2250,-1250){\makebox(0,0)[lb]{\smash{
{\color[rgb]{0,0,0}$\ba $}
}}}
\put(3300,-600){\makebox(0,0)[lb]{\smash{
{\color[rgb]{0,0,.69}$\bn_\theta $}
}}}
\put(3000,-500){\makebox(0,0)[lb]{\smash{
{\color[rgb]{1,0,0}$\bn $}
}}}
\put(2400,-800){\makebox(0,0)[lb]{\smash{
{\color[rgb]{0,0,.69}$\bt_\theta $}
}}}
\put(5200,-2300){\makebox(0,0)[lb]{\smash{
{\color[rgb]{0,0,0}${\cal B}_\theta^+ $}
}}}
\put(1500,-150){\makebox(0,0)[lb]{\smash{
{\color[rgb]{0,0,0}${\cal B}_\theta^- $}
}}}
\put(4100,-1100){\makebox(0,0)[lb]{\smash{
$\theta $
}}}
\end{picture}
}
\caption{If $\x_0 $ belongs to ${\cal B}_\theta^- $, we get mixed boundary conditions along $[\ba \bb ]$ and condition \eqref{S2} is satisfied.}
\end{figure}

\noindent We now can describe every situation by considering only three 
following cases (see Fig. 2),
$$
\hbox{(C1):} \quad 0<\theta_1 \le \theta_2 <\frac{\pi }{4} \, ,
\qquad
\hbox{(C2):} \quad 0<\theta_1 <\frac{\pi }{4} \le \theta_2 <\frac{\pi }{2} \, ,
\qquad
\hbox{(C3):} \quad \frac{\pi }{4}\le \theta_1 \le \theta_2 <\frac{\pi }{2} \, .
$$
\begin{figure}[htb]
\centerline{
\epsfig{figure=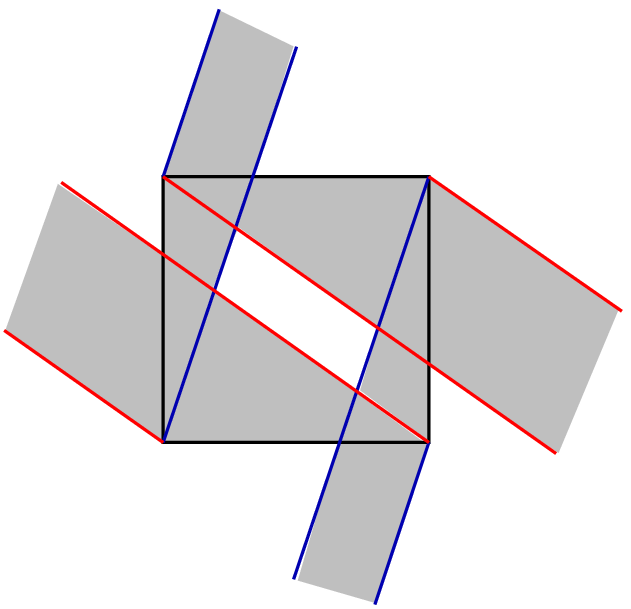,scale=0.70}
\qquad
\epsfig{figure=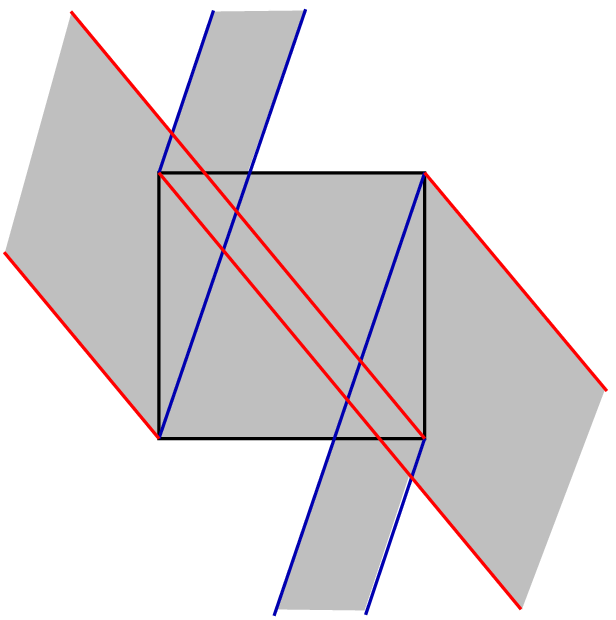,scale=0.70}
\qquad
\epsfig{figure=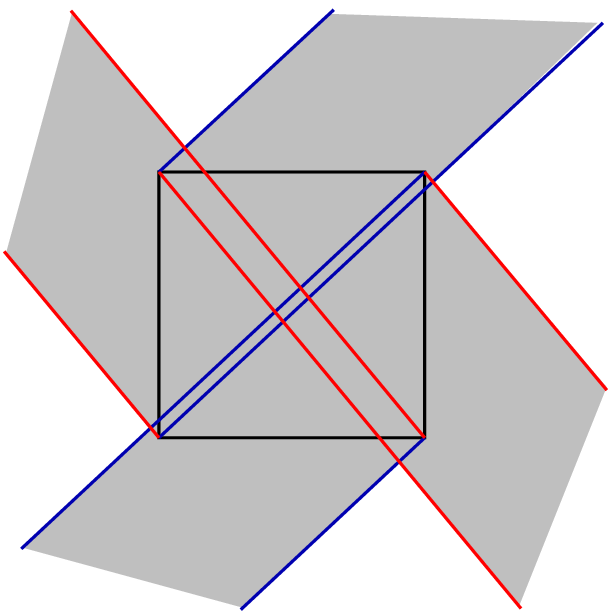,scale=0.70}
}
\caption{Cases (C1), (C2), (C3). Condition \eqref{S2} is not satisfied in
colored regions.}
\end{figure}

\noindent We also show a fully detailled partition in some particular
case coresponding to (C2) (see Fig. 3).

\begin{figure}[htb]
\centerline{
\begin{picture}(0,0)
\includegraphics{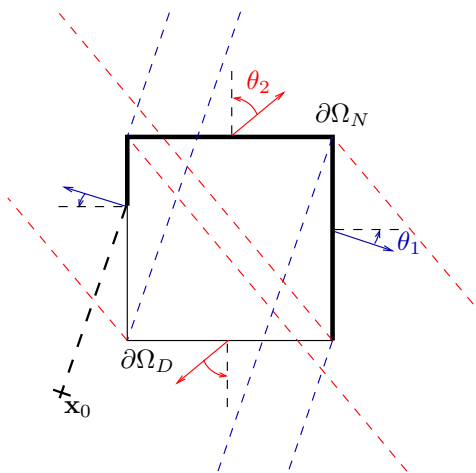}
\end{picture}
\setlength{\unitlength}{3729sp}
\begin{picture}(3085,3104)(1143,-3180)
\put(3540,-1670){\makebox(0,0)[lb]{\smash{
{\color[rgb]{0,0,.69}$\theta_1 $}
}}}
\put(2550,-600){\makebox(0,0)[lb]{\smash{
{\color[rgb]{1,0,0}$\theta_2 $}
}}}
\put(3000,-820){\makebox(0,0)[lb]{\smash{
{\color[rgb]{0,0,0}$\partial \Omega_N $}
}}}
\put(1720,-2450){\makebox(0,0)[lb]{\smash{
{\color[rgb]{0,0,0}$\partial \Omega_D $}
}}}
\put(1350,-2750){\makebox(0,0)[lb]{\smash{
{\color[rgb]{0,0,0}$\x_0 $}
}}}
\end{picture}
}
\caption{Example of Dirichlet and Neumann parts of the boundary in a case of
(C2)-type.}
\end{figure}

\end{document}